\documentclass[12pt]{article}

\usepackage{graphicx}
\usepackage{amsmath}
\usepackage{amsthm}
\usepackage{amssymb}
\usepackage[mathscr]{eucal}
\usepackage[all]{xy}
\usepackage{exscale,relsize}
\usepackage{psfrag}

\title{\bf A computation in Khovanov-Rozansky Homology}
\author{Daniel Krasner}
\date{}

\theoremstyle{plain}
\newtheorem{theorem}{Theorem}

\newtheorem{lemma}[theorem]{Lemma}

\newtheorem{question}[theorem]{Question}

\def\title{\em}

\usepackage[top=1in, bottom=1in, left=1in, right=1in]{geometry}

\begin{document}

\newpage

\maketitle

\medskip
\noindent

\begin {abstract}
We investigate the Khovanov-Rozansky invariant of a certain tangle and its compositions. Surprisingly the complexes we
encounter reduce to ones that are very simple. Furthermore, we
discuss a ``local" algorithm for computing Khovanov-Rozansky homology
and compare our results with those for the ``foam" version of 
$sl_3$-homology.
\end{abstract}
\newpage
\medskip\noindent

\section{Introduction}

\begin{figure}[h]
\centerline{
\includegraphics[scale = .8]{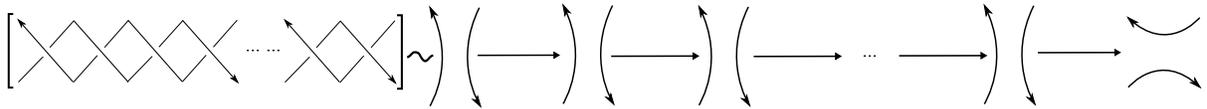}}
\caption{Our main tangle and its reduced complex} 
\label{intro}
\end{figure}

In a seminal work M. Khovanov and L. Rozansky \cite{KR} introduced a
series of doubly-graded link homology theories with Euler
characteristic the quantum $sl_{n}$-link polynomials. The
construction relied on the theory of matrix factorizations, which
was previously seen in the study of maximal Cohen-Macaulay modules
on isolated hypersurface singularities. For $n=2$ and $n=3$,  link
homology theories with Euler characteristic the Jones polynomial and
the quantum $sl_{3}$ polynomial respectively, were introduced
earlier by M. Khovanov in \cite{Kh2} and \cite{Kh1}. The
constructions came in a very different guise, but it was easy to see
that the matrix factorization version specialized to $n=2$ agreed
with what is now known as Khovanov homology. The $sl_{3}$ version is
also know to be isomorphic to the the matrix factorization version \cite{VM2} . Variants of these theories were described in \cite{BN1}, \cite{BN2},
\cite{VM} as well as a number of other publications. Using ideas
from \cite{BN3} we show that for certain classes of tangles, and hence for knots and links composed of these, the
Khovanov-Rozansky complex reduces to one that is quite simple,
that is one without any ``thick" edges. In particular we consider the tangle in figure \ref{intro} and show that its associated complex is homotopic to the one below, with some grading shifts and basic maps which we leave out for now.

The complexes for these knots and links are entirely ``local," and to calculate the homology  we only need to exploit the Frobenius structure of the underlying algebra assigned to the unknot. Hence, here the calculations and complexity is similar to that of $sl_2$-homology. We also discuss a general
algorithm, basically the one described in \cite{BN3}, to compute
these homology groups in a more time-efficient manner. We compare our
results with similar computations in the version of $sl_{3}$-homology found in \cite{Kh1}, which we refer to as the "foam"
version (foams are certain types of cobordisms described in this
paper), and giving an explicit isomorphism between the two versions. A very similar calculation in the $sl_{3}$-homology, that for the $(2,n)$ torus knots, was first done in \cite{MN}. The paper
is structured as follows: in section $2$ we give a brief review of
Khovanov-Rozansky homology, but assume that the reader is either
familiar with the material or is willing to take a lot for granted;
in section $3$ through $5$ we go through the main calculation; in
section $6$ we discuss the algorithm and "foam" version of 
$sl_{3}$-homology.

\medskip

\textbf{Acknowledgements}

Firstly, and above all, I would like to thank my advisor Mikhail Khovanov.
 I would also like to acknowledge Yanfeng Chen for helpful discussion and Scott Morrison for pointing me to Bar-Natan's paper \cite{BN2} and to the calculations in \cite{MN}. In addition many thanks to Jacob Rasmussen for his many helpful suggestions on the first draft.  

\pagebreak

\section{A Review of Khovanov-Rozansky Homology}\

\medskip

\textbf{Matrix Factorizations}\

Let $R =  \mathbb{Q}[x_{1}, \dots , x_{n}]$ be a graded polynomial ring
 in $n$ variables with $\deg(x_i)=2$, and let $\omega \in R$. A \emph{matrix
factorization} with \emph{potential} $\omega$ is a collection of two
free $R$-modules $M^{0}$ and $M^{1}$ and $R$-module maps $d^{0}: M^{0}
\rightarrow M^{1}$ and $d^{1}: M^{1} \rightarrow M^{0}$ such that:

\begin{center}
$d^{0} \circ d^{1} = \omega$ Id and $d^{1} \circ d^{0} =
\omega$ Id
\end{center}

The $d^{i}$'s are referred to as 'differentials' and we often denote
such a 2-complex by

\begin{displaymath}
\xymatrix{
 M = & M^{0}   \ar[r]^{d^{0}}   &  M^{1}  \ar[r]^{d^{1}} & M^{0}}
\end{displaymath}

Given two matrix factorizations $M_{1}$ and $M_{2}$ with potentials
$\omega_{1}$ and $\omega_{2}$ respectively,  their tensor product
 is given as the tensor product of complexes, and it is easy to see
 that $M_{1} \otimes M_{2}$ is a matrix factorization with
 potential $\omega_{1} + \omega_{2}$.

To keep track of minus signs, it is convenient
to assign a label to the factorization and denote it by

\begin{displaymath}
\xymatrix{
 M = & M(\varnothing)   \ar[r]^{d^{0}}   &  M(a)  \ar[r]^{d^{1}} & M(\varnothing),}
\end{displaymath}

so that the tensor product of two factorizations $M \otimes M$ can be written as

\begin{center}
$\left(
\begin{array}{clcr}
M(\varnothing)\\
M(ab)
\end{array} \right)$
$\stackrel{}{\longrightarrow}$ $\left(
\begin{array}{clcr}
M(a)\\
M(b)
\end{array} \right)$
$\stackrel{}{\longrightarrow}$ $\left(
\begin{array}{clcr}
M(\varnothing)\\
M(ab)
\end{array} \right).$
\end{center}

Here we are simply replacing $M^0$ by $M(\varnothing)$ and $M^1$ by a label such as $M(a)$; this will be useful below when we assign facorizations to plane graphs. See \cite{KR} for a more detailed treatement.  

 A homomorphism $f: M \rightarrow N$ of two factorizations is a pair
 of homomorphisms $f^{0}: M^{0} \rightarrow N^{0}$ and $f^{1}: M^{1} \rightarrow
 N^{1}$ such that the following diagram is commutative:

\begin{displaymath}
\xymatrix{
M^{0}  \ar[d]^{f^{0}}   \ar[r]^{d^{0}} & M^{1}  \ar[d]^{f^{1}} \ar[r]^{d^{1}} &  M^{0} \ar[d]^{f^{0}}\\
N^{0}   \ar[r]^{d^{0}}   &  N^{1}  \ar[r]^{d^{1}} & N^{0}}
\end{displaymath}

A homotopy $h$ between maps $f, g: M \rightarrow N$ of factorizations
is a pair of maps $h^{i}: M^{i} \rightarrow N^{i-1}$ such that $f -
g = h \circ d_{M} + d_{N} \circ h$ where $d_{M}$ and $d_{N}$ are the
differentials in $M$ and $N$ respectively. For a detailed treatment
of matrix factorizations we refer the reader to \cite{KR}.
\medskip

\textbf{Grading Shifts}

Let $M$ be a matrix factorization as above, with $M^0$ and $M^1$ $\mathbb{Z}$-graded modules over a $\mathbb{Z}$-graded
ring and let $k \in \mathbb{Z}$. Let $M\{k\}$ be the module $M$
with degrees shifted up by $k$.
By $M\langle k \rangle^{i}  = M^{i+k}$ with $i+k$ taken mod 2 we denote the shift in homological grading coming from the factorization. Later we will see another homological grading of our complex, arising from the resolutions of a
link diagram, and the shifted module there will be denoted by $M[k]$.
\medskip

\textbf{Planar Graphs and Matrix Factorizations}

\begin{figure}[h]
\centerline{
\includegraphics[scale=.7]{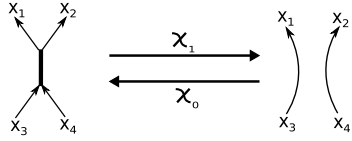}}
\caption{Maps $\chi_{0}$ and $\chi_{1}$}
\label{maps}
\end{figure}

Our graphs are embedded in a disk and have two types of edges,
unoriented and oriented. Unoriented edges are called ``thick"
and drawn accordingly; each vertex adjoining a thick edge has either
two oriented edges leaving it or two entering. In figure \ref{maps} left $x_1, x_2$ are outgoing and $x_3, x_4$ are incoming. Oriented edges are
allowed to have marks and we also allow closed loops; points of the
boundary are also referred to as marks. See for example figure \ref{A planar graph}
below. To such a graph $\Gamma$ we assign a matrix factorization in
the following manner:

To a thick edge $t$ as in figure \ref{maps} left we assign a factorization $C_{t}$
with potential $\omega_{t} = x_{1}^{n+1} + x_{2}^{n+1} - x_{3}^{n+1}
- x_{4}^{n+1}$ over the ring $R_{t}= \mathbb{Q}[x_{1}, x_{2},
x_{3}, x_{4}]$. Since $x^{n+1} + y^{n+1}$ lies in the ideal
generated by $x+y$ and $xy$ we can write it as a polynomial $g(x+y,
xy)$. Hence, $\omega_{t}$ can be written as

\begin{center}
$\omega_{t} = (x_{1} + x_{2} - x_{3} - x_{4})u_{1} +
(x_{1}x_{2}-x_{3}x_{4})u_{2}$
\end{center}

where

$$u_{1} = \displaystyle \frac{x_{1}^{n+1}+ x_{2}^{n+1} - g(x_{3} +
x_{4}, x_{1}x_{2})}
{x_{1} + x_{2} - x_{3} - x_{4}},$$ 

$$u_{2} = \displaystyle \frac{g(x_{3} + x_{4}, x_{1}x_{2}) -
x_{3}^{n+1} - x_{4}^{n+1}}{ x_{1}x_{2}-x_{3}x_{4}}.$$

$C_{t}$ is the tensor product of graded
factorizations

\begin{center}
$R_{t} \xrightarrow{u_{1}}
    R_{t}\{1-n\} \xrightarrow{x_{1}+ x_{2} - x_{3} - x_{4}}  R_{t}$
\end{center}
and

$$R_{t} \xrightarrow{u_{2}}  R_{t}\{3-n\} \xrightarrow{x_{1}x_{2} -
x_{3}x_{4}} R_{t}.$$

To an arc $\alpha$ bounded by marks oriented from $j$ to $i$ we
assign the factorization $L_{j}^{i}$

$$R_{\alpha}   \xrightarrow{\pi_{ij}}
    R_{\alpha} \xrightarrow{x_{i} - x_{j}}  R_{\alpha}, $$

where $R_{\alpha} = \mathbb{Q}[x_{i}, x_{j}]$ and

$$\pi_{ij} = \displaystyle \frac{x_{i}^{n+1} - x_{j}^{n+1}}{{x_{i} -
x_{j}}}.$$

Finally, to an oriented loop with no marks we assign the complex $0
\rightarrow A \rightarrow 0 = A\langle 1 \rangle$ where $A =
\mathbb{Q}[x]/(x^{n})$. [Note: to a loop with marks we 
assign the tensor product of $L_{j}^{i}$'s as above, but this
turns out to be isomorphic  to $A\langle 1 \rangle$ in the homotopy category.]

\begin{figure}[h]
\centerline{
\includegraphics[scale=.7]{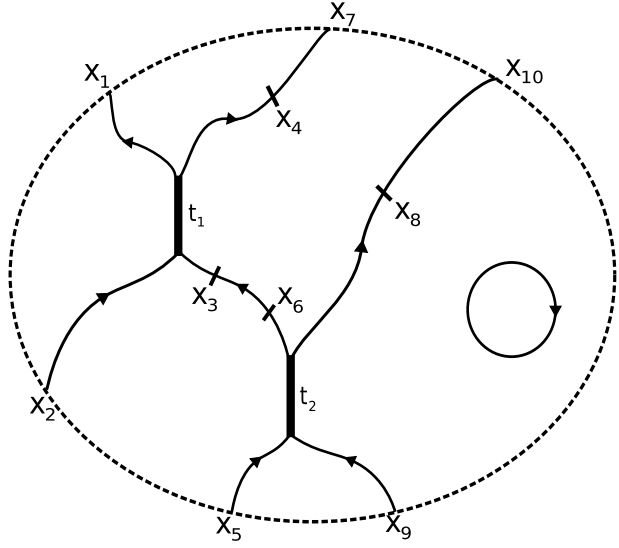}}
\caption{A planar graph}
\label{A planar graph}
\end{figure}

We define $C(\Gamma)$ to be the tensor product of $C_{t}$ over all
thick edges $t$, $L^{i}_{j}$ over all edges $\alpha$ from $i$ to $j$,
and $A\langle 1 \rangle$ over all oriented markless loops.
This tensor product is taken over appropriate rings such that
$C[\Gamma]$ is a free module over $R= \mathbb{Q}[x_{i}]$ where the
$x_{i}$'s are marks. For example to the graph in
figure \ref{A planar graph} we assign $C(\Gamma) =
L^{7}_{4}  \otimes C_{t_{1}} \otimes L^{3}_{6}  \otimes C_{t_{2}} \otimes L^{10}_{8} \otimes A\langle 1
\rangle$ tensored over $\mathbb{Q}[x_{4}]$, $\mathbb{Q}[x_{3}]$,
$\mathbb{Q}[x_{6}]$, $\mathbb{Q}[x_{8}]$ respectively.
$C(\Gamma)$ becomes a $\mathbb{Z} \oplus
\mathbb{Z}_{2}$-graded complex with the $\mathbb{Z}_{2}$-grading
coming from the factorization. It has potential $\omega =
\displaystyle \sum_{i \in \partial \Gamma} \pm x_{i}^{n+1}$, where
$\partial \Gamma$ is the set of all boundary marks and the $+$, $-$
is determined by whether the direction of the edge corresponding to
$x_{i}$ is towards or away from the boundary. [Note: if $\Gamma$ is
a closed graph the potential is zero.]
\medskip

\textbf{The maps $\chi_{0}$ and $\chi_{1}$}

We now define maps between matrix factorizations associated to
the thick edge and two disjoint arcs as in figure \ref{maps}. Let
$\Gamma^{0}$ correspond to the two disjoint arcs and $\Gamma^{1}$ to
the thick edge.

$C(\Gamma^{0})$ is the tensor product of $L_{4}^{1}$ and
$L_{3}^{2}$. If we assign labels $a$, $b$ to $L_{4}^{1}$,
$L_{3}^{2}$ respectively, the tensor product can be written as

$$\left(
\begin{array}{clcr}
R(\varnothing)\\
R(ab)\{2-2n\}
\end{array} \right)
\stackrel{P_{0}}{\longrightarrow}
\left(
\begin{array}{clcr}
R(a)\{1-n\}\\
R(b)\{1-n\}
\end{array} \right)
\stackrel{P_{1}}{\longrightarrow}
\left(
\begin{array}{clcr}
R(\varnothing)\\
R(ab)\{2-2n\}
\end{array} \right),$$

where

$$P_{0} = \left(
\begin{array}{cclcr}
\pi_{14} & x_{2} - x_{3}\\
\pi_{23} & x_{4} - x_{1}
\end{array} \right)
, P_{1} = \left(
\begin{array}{cclcr}
x_{1} - x_{4} & x_{2} - x_{3}\\
\pi_{23} & -\pi_{14}
\end{array} \right),$$

$$\pi_{ij} = \displaystyle\sum_{k=0}^{n} x_{i}^{k}x_{j}^{n-k}.$$

Assigning labels $a'$ and $b'$ to the two factorizations in 
$C(\Gamma^{1})$, we have that $C(\Gamma^{1})$ is given by

$$\left(
\begin{array}{clcr}
R(\varnothing)\{-1\}\\
R(a'b')\{3-2n\}
\end{array} \right)
\stackrel{Q_{1}}{\longrightarrow} 
\left(
\begin{array}{clcr}
R(a')\{n\}\\
R(b')\{2-n\}
\end{array} \right)
\stackrel{Q_{2}}{\longrightarrow} 
\left(
\begin{array}{clcr}
R(\varnothing)\{-1\}\\
R(a'b')\{3-2n\}
\end{array} \right),$$

where

$$Q_{1} = \left(
\begin{array}{cclcr}
u_{1} & x_{1}x_{2} - x_{3}x_{4}\\
u_{2} & x_{3} + x_{4} - x_{1} - x_{2}
\end{array} \right)
, Q_{2} = \left(
\begin{array}{cclcr}
x_{1} + x_{2} - x_{3} - x_{4} & x_{1}x_{2} - x_{3}x_{4}\\
u_{2} & -u_{1}
\end{array} \right).$$

A map between $C(\Gamma^{0})$ and $C(\Gamma^{1})$ can be given by a
pair of $2\times2$ matrices. Define $\chi_{0}: C(\Gamma^{0})
\rightarrow C(\Gamma^{1})$ by

$$U_{0} = \left(
\begin{array}{cclcr}
x_{1} - x_{3} & 0\\
\frac{u_{1}+x_{1}u_{2}-\pi_{23}}{x_{1}-x_{4}} & 1
\end{array} \right)
, U_{1} = \left(
\begin{array}{cclcr}
x_{1} & -x_{3}\\
-1 & 1
\end{array} \right),$$

and $\chi_{1}: C(\Gamma^{1}) \rightarrow C(\Gamma^{0})$  by

$$V_{0} = \left(
\begin{array}{cclcr}
1 & 0\\
\frac{u_{1}+x_{1}u_{2}-\pi_{23}}{x_{4}-x_{1}} & x_1 - x_3
\end{array} \right)
, V_{1} = \left(
\begin{array}{cclcr}
1 & x_{3}\\
1 & x_{1}
\end{array} \right).$$

These maps have degree $1$. 
Computing we see that the composition  $\chi_{1}\chi_{0} = (x_{1}- x_{3})I$, where $I$ is the identity matrix, i.e.
$\chi_{1}\chi_{0}$ is multiplication by $x_{1} - x_{3}$. Similarly
$\chi_{0}\chi_{1} = (x_{4} - x_{2})I$. [Note: these are
specializations of the  maps $\chi_{0}$ and $\chi_{1}$ given in
\cite{KR}, with $\lambda = 0$ and $\mu = 1$. As these maps
are homotopic for any rational value of $\lambda$ and $\mu$ we are free to
do so.]\\

 Define the trace $\varepsilon:\mathbb{Q}[x]/(x^n) \longrightarrow \mathbb{Q}$ as 
$\varepsilon(x^i)=0$ for $i \neq n-1$ and $\varepsilon(x^{n-1})=1$. The 
unit $\iota: \mathbb{Q}  \longrightarrow \mathbb{Q}[x]/(x^n)$ is defined 
by $\iota(1)=1$.\\

The relations between $C(\Gamma)$'s  mimic the graph skein relations, see for example
\cite{KR}, and we list the ones needed below.

\textbf{Direct Sum Decomposition 0:}

\medskip
\begin{figure}[h]
\centerline{
\includegraphics[scale=.9]{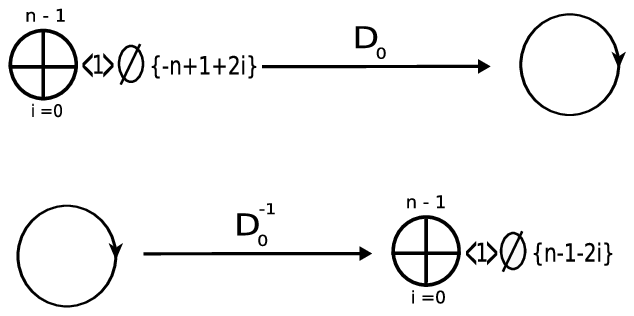}}
\end{figure}

where $D_{0} = \displaystyle\sum_{i=0}^{n-1} x^{i} \iota$ and
$D_{0}^{-1} = \displaystyle\sum_{i=0}^{n-1} \varepsilon x^{n-1-i}.$

By the pictures above, we really mean the complexes assigned to them, i.e. $\emptyset \langle 1 \rangle$ is the complex with $\mathbb{Q}$ sitting in homological grading $1$ and the unknot is the complex $A \langle 1 \rangle$ as above. The map $x^i \iota$ is a composition of maps 

$$A \langle 1 \rangle \xrightarrow{x^i} \langle 1 \rangle \xrightarrow{\iota} \emptyset \langle 1 \rangle,$$

where $x^i$ is multiplication and $\iota$ is the unit map, i.e. $x^i \iota$ is the map

$$\mathbb{Q}[x]/(x^n) \xrightarrow{x^i} \mathbb{Q}[x]/(x^n) \xrightarrow{\iota} \mathbb{Q} .$$ 

Similar with $ \varepsilon x^{n-1-i}$.
It is easy to check that the above maps are grading preserving and
their composition is the identity.

\bigskip

 \textbf{Direct Sum Decomposition I:}\
\medskip

\begin{figure}[h]
\centerline{
\includegraphics[scale=.9]{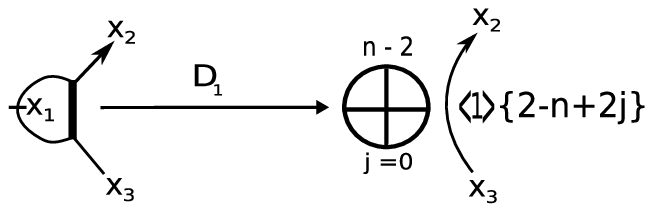}}
\end{figure}
\begin{figure}[h]
\centerline{
\includegraphics[scale=.8]{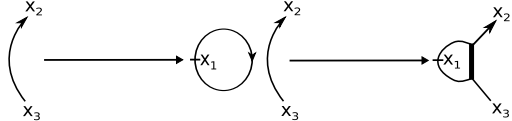}}
\caption{The map $\alpha$}
\label{alpha}
\end{figure}
\begin{figure}[h]
\centerline{
\includegraphics[scale=.8]{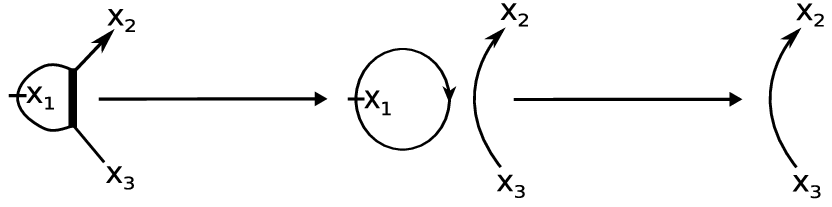}}
\caption{The map $\beta$}
\label{beta}
\end{figure}


where $D_{1} = \displaystyle\sum_{i=0}^{n-2} \beta x_{1}^{n-i-2}$ and $D_{1}^{-1} = \displaystyle\sum_{i=0}^{n-2}
\sum_{j=0}^{i} x_{1}^{j}x_{2}^{i-j} \alpha$
with $\alpha:= \chi_{0}  \circ \iota'$ and $\beta:= \varepsilon' \circ  \chi_{1}$. Here $\iota' = \iota \otimes Id$ and $\varepsilon' = \varepsilon \otimes Id$; the $Id$ corresponds to the arc with endpoints labeled by $x_2, x_3$, i.e  $\iota'$ is the map that includes the single arc diagram into one with the unkot and single arc disjoint, see figure \ref{alpha}. Similar with $\varepsilon'$ in the right half of figure \ref{beta}.

\medskip

\medskip

\pagebreak

\textbf{Direct Sum Decomposition II:}\
\medskip

\begin{figure}[h]
\centerline{
\includegraphics[scale=.75]{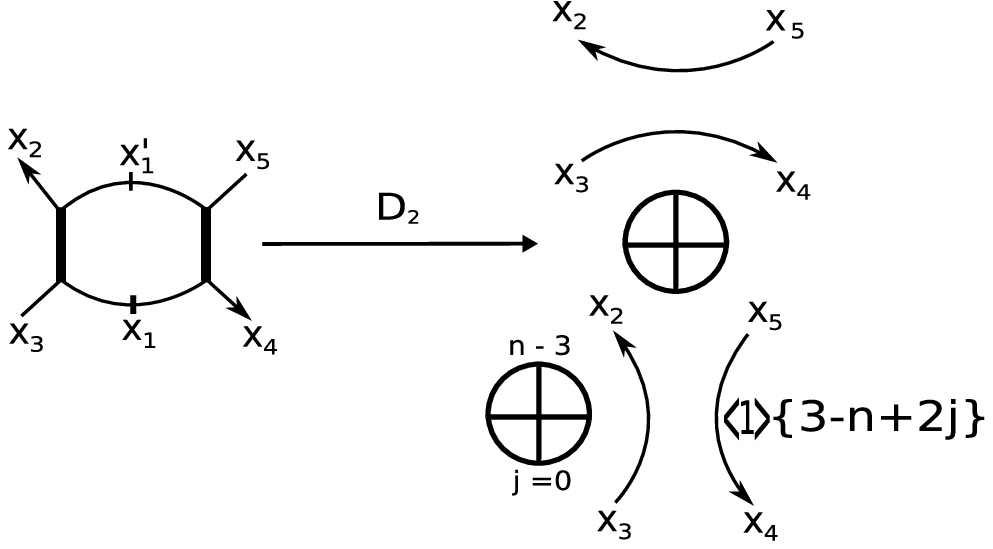}}
\end{figure}

where $D_{2} = S \oplus
 \displaystyle\sum_{j=0}^{n-3}\beta_{j}$
 and $\beta_{j}= \displaystyle\sum_{j=0}^{n-3}\beta
\sum_{a+b+c=n-3-j}x_{2}^{a}x_{4}^{b}x_{1}^{c}.$

\begin{figure}[h]
\centerline{
\includegraphics[scale=.7]{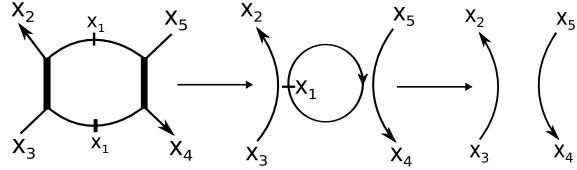}}
\caption{The map $\beta$}
\label{beta2}
\end{figure}

Here $\beta$ is given by the composition of two $\chi_1$'s, corresponding to the two thick edges on the left-hand side above, and the trace map $\varepsilon$, see figure \ref{beta2}. Finally $S$ is gotten by ``merging" the thick edges together to form two disjoint horizontal arcs, as in the top righ-hand corner above; an exact descrition of $S$ won't really matter so we will not go into details and refer the interested reader to \cite{KR}. 

\pagebreak

\textbf{Tangles and complexes}

\medskip

We resolve a crossing $p$ in the two ways and assign to it a complex
$C^{p}$ depending on whether the crossing is positive or
negative. To a diagram $D$ representing a tangle $L$ we assign the
complex $C(D)$ of matrix factorization which is the tensor product
of $C^{p}$, over all crossings $p$, of $L_{j}^{i}$ over arcs $j
\rightarrow i$, and of $A \langle 1\rangle$ over all crossingless
markless circles in $D$. The tensor product is taken as before so
that $C(D)$ is free and of finite rank as an $R$-module. This
complex is $\mathbb{Z} \oplus \mathbb{Z} \oplus \mathbb{Z}_{2}$
graded.

\begin{figure}[h]
\centerline{
\includegraphics[scale=.7]{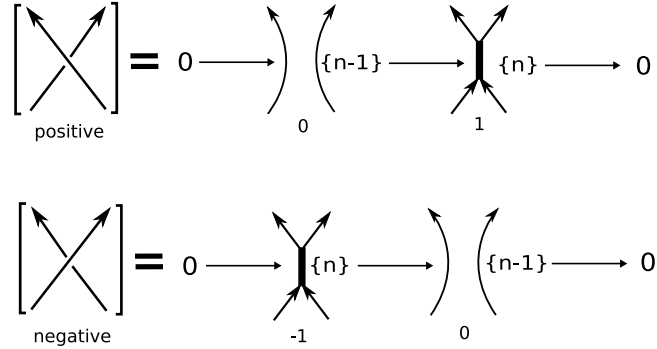}}
\caption{Complexes associated to pos/neg crossings; the
numbers below the diagrams are cohomological degrees.}
\label{crossings}
\end{figure}

\begin{theorem}
(Khovanov-Rozansky, \cite{KR})
 The isomorphism class of
$C(D)$ up to homotopy is an invariant of the tangle.
\end{theorem}

If $L$ is a link the cohomology groups are nontrivial only in degree
equal to the number of components of $L$ mod $2$. Hence, the grading
reduces to $\mathbb{Z} \oplus \mathbb{Z}$. The resulting cohomology
groups are denoted by

\begin{center}
$H_{n}(D) = \displaystyle \bigoplus_{i,j \in \mathbb{Z}}H^{i,j}_{n}(D),$
\end{center}

and the Euler characteristic of $H_{n}(D)$ is the quantum link
polynomial $P_{n}(L)$, i.e.

\begin{center}
$P_{n}(L) = \displaystyle \sum_{i,j \in \mathbb{Z}} (-1)^{i}q^{j}
dim_{\mathbb{Q}} H^{i,j}_{n}(D).$
\end{center}

The isomorphism classes of $H^{i,j}_{n}(D)$ depend only on the link
$L$ and, hence, are invariants of the link.
\medskip

\textbf{Gaussian Elimination for Complexes:}\

\begin{lemma}  If $\phi:B \rightarrow D$ is an isomorphism (in some
additive category $\cal C$), then the four term complex segment
below

\begin{equation}
  \xymatrix@C=2cm{
    \cdots\
    \left[A\right]
    \ar[r]^{\begin{pmatrix}\alpha \\ \beta\end{pmatrix}} &
    {\begin{bmatrix}B \\ C\end{bmatrix}}
    \ar[r]^{\begin{pmatrix}
      \phi & \delta \\ \gamma & \epsilon
    \end{pmatrix}} &
    {\begin{bmatrix}D \\ E\end{bmatrix}}
    \ar[r]^{\begin{pmatrix} \mu & \nu \end{pmatrix}} &
    \left[F\right] \  \cdots
  }
\end{equation}
is isomorphic to the (direct sum) complex segment
\begin{equation}
  \xymatrix@C=3cm{
    \cdots\
    \left[A\right]
    \ar[r]^{\begin{pmatrix}0 \\ \beta\end{pmatrix}} &
    {\begin{bmatrix}B \\ C\end{bmatrix}}
    \ar[r]^{\begin{pmatrix}
      \phi & 0 \\ 0 & \epsilon-\gamma\phi^{-1}\delta
    \end{pmatrix}} &
    {\begin{bmatrix}D \\ E\end{bmatrix}}
    \ar[r]^{\begin{pmatrix} 0 & \nu \end{pmatrix}} &
    \left[F\right] \  \cdots
  }.
\end{equation}
Both of these complexes are homotopy equivalent to the (simpler)
complex segment
\begin{equation}
  \xymatrix@C=3cm{
    \cdots\
    \left[A \right]
    \ar[r]^{\left(\beta\right)} &
    {\left[C \right]}
    \ar[r]^{\left(\epsilon-\gamma\phi^{-1}\delta\right)} &
    {\left[E\right]}
    \ar[r]^{\left(\nu\right)} &
    \left[F\right] \  \cdots
  }.
\end{equation}
Here the capital letters are arbitrary columns of objects in $\cal C$ and all
Greek letters are arbitrary matrices representing morphisms with the
appropriate dimensions, domains and ranges (all the matrices are
block matrices); $\phi : B \rightarrow D$ is an isomorphism, i.e. it
is invertible.
\end{lemma}

\emph{Proof:} The matrices in complexes $(1)$ and $(2)$ differ by a
change of bases, and hence the complexes are isomorphic. $(2)$ and
$(3)$ differe by the removal of a contractible summand; hence, they
are homotopy equivalent. $\square$

\pagebreak

\section{The Basic Calculation}

\begin{figure}[h]
\centerline{
\includegraphics[scale=.75]{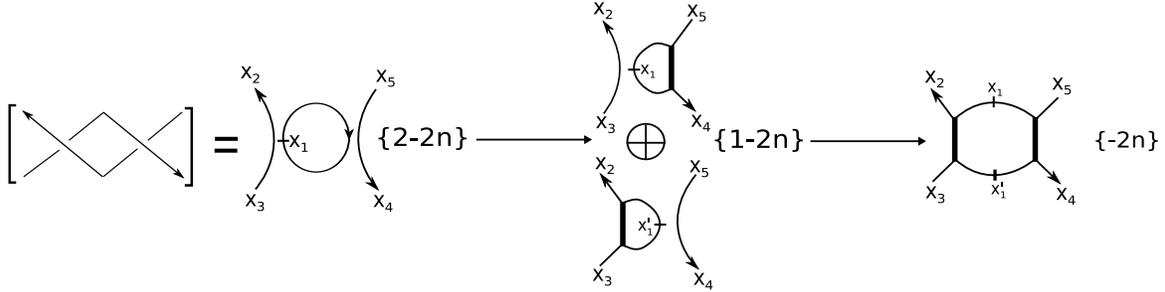}}
\caption{The tangle $T$ and its complex}
\label{basictangle}
\end{figure}

We first consider the complex associated to the tangle $T$ in figure \ref{basictangle} with the appropriate maps $\chi_0$ and $\chi_1$ left out.

\ We first look at the following part of the complex and, for
simplicity, leave out the overall grading shifts until later:

\begin{figure}[h]
\centerline{
\includegraphics[scale=.8]{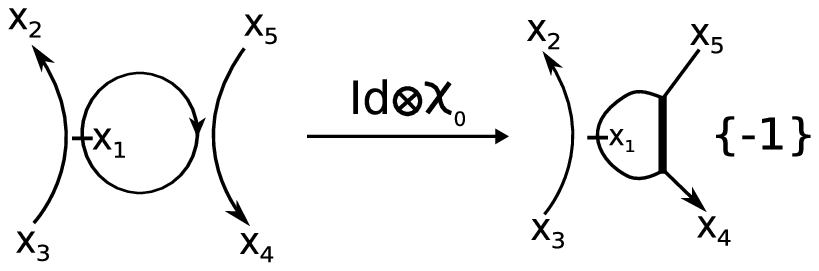}}
\end{figure}

We apply direct sum decompositions 0 and I and end up with the
following where the maps $F_{1}$ and $F_{2}$ are isomorphisms:

\begin{figure}[h]
\centerline{
\includegraphics[scale=.9]{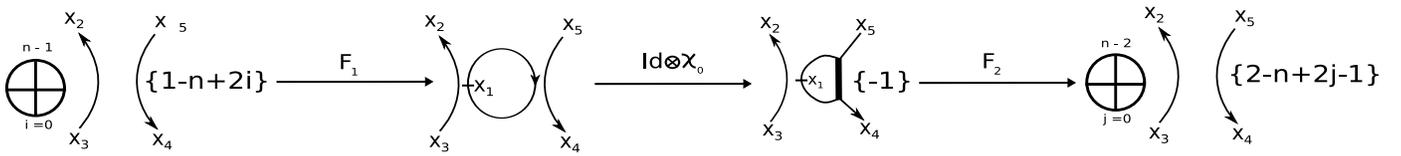}}
\caption{First part of the complex for T with decompositions}
\label{Tfig2}
\end{figure}

Explicitly,
 $F_{1} = \sum_{i=0}^{n-1} Id \otimes x^{i}_{1} \iota \otimes Id$
and $F_{2} = \sum_{j=0}^{n-2} Id \otimes \beta_{j}$ \

Composing the maps we get:
\begin{small}
$$\begin{array}{rcl} F_{2} \circ (Id \otimes \chi_{0}) \circ F_{1} &
= & (\displaystyle\sum_{j=0}^{n-2} Id \otimes \beta_{j}) \circ (Id
\otimes \chi_{0}) \circ (\sum_{i=0}^{n-1} Id \otimes x_{1}^{i} \iota
\otimes
Id)\\
& = & (\displaystyle\sum_{j=0}^{n-2} Id \otimes \beta_{j}) \circ
(\sum_{i=0}^{n-1} Id
\otimes (\chi_{0} \circ (x_{1}^{i} \iota \otimes Id)))\\
 & = & \displaystyle\sum_{j=0}^{n-2} \sum_{i=0}^{n-1} Id \otimes (\beta_{j} \circ
\chi_{0} \circ (x_{1}^{i} \iota \otimes Id))\\
& = & \displaystyle\sum_{j=0}^{n-2} \sum_{i=0}^{n-1} Id \otimes
(\varepsilon'(x_{1} -
x_{4})x_{1}^{n+i-j-2})\\
 & = & \displaystyle\sum_{j=0}^{n-2} \sum_{i=0}^{n-1} Id \otimes
(\varepsilon'(x_{1}^{n+i-j-1} - x_{4}x_{1}^{n+i-j-2}))\\
& = & \displaystyle\sum_{j=0}^{n-2} \sum_{i=0}^{n-1} Id \otimes
\underbrace{(\varepsilon(x_{1}^{n+i-j-1}) -
x_{4}\varepsilon(x_{1}^{n+i-j-2}))}_{\Theta}.
\end{array}$$
\end{small}

To go from line 3 to 4 and 4 to 5, recall that $\beta_{j}=
\varepsilon' \circ \chi_{1} x_{1}^{n-j-2}$ and $\chi_{1} \circ
\chi_{0} = x_{1} - x_{4} = x_{1} - x_{5}$. [Note: for lack of better notation, we use ``$\sum$" to indicate both a map from a direct sum and an actual sum, as seen above indexed $i$ and $j$ respectively.] 

 Now $\Theta = Id$ if $i=j$, $-x_{4}$ if $i=j+1$, and $0$ otherwise,
  $F_{2} \circ (Id \otimes \chi_{0}) \circ F_{1}$ is given by the
 following $(n-1) \times n-1$ matrix:
 \\
\begin{small}
\centerline {$\left[
\begin{array}{ccclcr}
Id & -x_{4} & 0 &  \ldots & \ldots &  0\\
0  &  Id    & -x_{4} & 0 & \ldots & 0 \\
\vdots & & \ddots & \ddots & & \vdots\\
\vdots & &  & \ddots & \ddots& \vdots\\
0 & & \ldots & 0 & Id & -x_{4}\\
\end{array} \right]$}
\end{small}
 \
 \

\indent Using Gaussian Elimination for complexes it is easy to see
that, up to homotopy, only the top degree term survives. By degree,
we mean with respect to the above grading shifts.

\

\indent Now we look at the following subcomplex:

\begin{figure}[h]
\centerline{
\includegraphics[scale=.35]{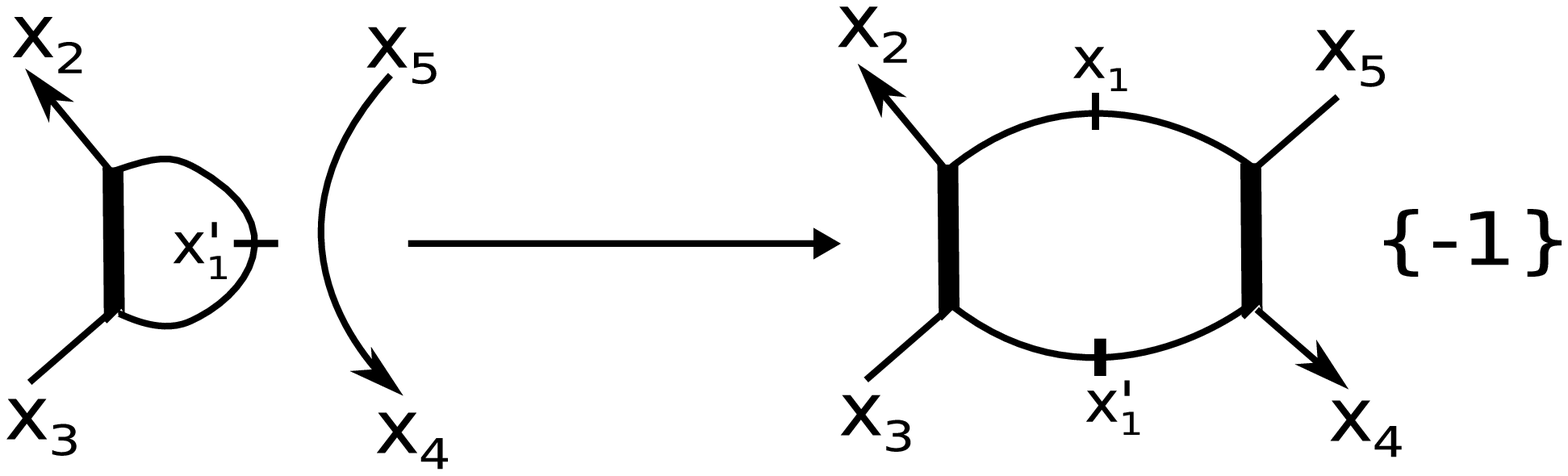}}
\end{figure}

\begin{figure}[h]
\centerline{
\includegraphics[scale=.8]{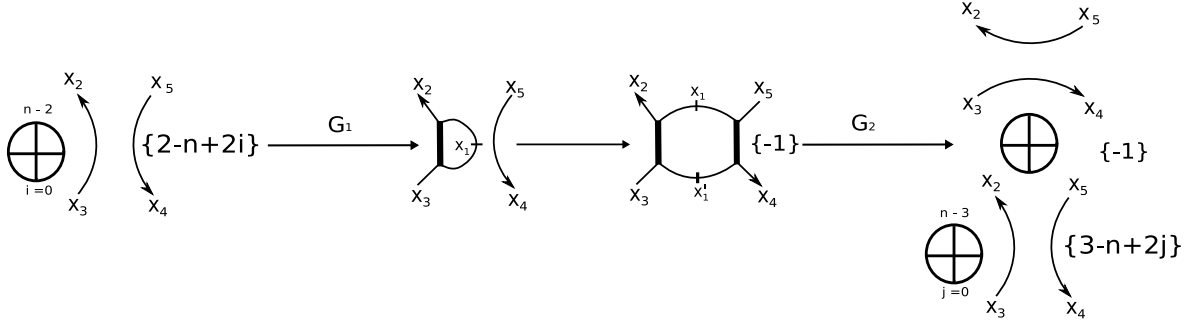}}
\caption{The second part of the complex for T with decompositions}
\label{Tfig4}
\end{figure}

 Including all the isomorphisms we have the complex in figure \ref{Tfig4},
with $G_{1} = \sum_{i=0}^{n-2} \alpha_{i} \otimes Id$ and $G_{2} =
S \oplus \sum_{j=0}^{n-3}\beta_{j}$ ($S$ is the saddle map).\

Composing these maps we get:

$$
\begin{array}{rcl} G_{2} \circ \chi_{0}^{''} \circ G_{1}  & 
= & (S \oplus \displaystyle\sum_{j=0}^{n-3}\beta_{j}) \circ \chi_{0}^{''} \circ( \sum_{i=0}^{n-2} \alpha_{i} \otimes Id)\\
 
& = & \left(S \oplus  \displaystyle\sum_{j=0}^{n-3}\beta
\sum_{a+b+c=n-3-j}x_{2}^{a}x_{4}^{b}x_{1}^{c}\right) \circ \chi_{0}^{''} \circ \left( \displaystyle\sum_{i=0}^{n-2} \displaystyle\sum_{k=0}^{i} x_{1}^{k}x_2^{i-k} \alpha \otimes Id \right)\\

& = & \left(S \oplus  \displaystyle\sum_{j=0}^{n-3}\varepsilon'
\circ \chi_{1}^{''} \circ \chi_{1}^{'} \sum_{a+b+c=n-3-j}x_{2}^{a}x_{4}^{b}x_{1}^{c}\right) \circ
\chi_{0}^{''} \circ \left( \displaystyle\sum_{i=0}^{n-2} \displaystyle\sum_{k=0}^{i} x_{1}^{k}x_2^{i-k}
\chi_{0}^{'} \circ \iota' \otimes Id \right)\\

& = & \overline{S} \oplus  \displaystyle\sum_{j=0}^{n-3} \sum_{i=0}^{n-2}
\varepsilon' \circ \chi_{1}^{''} \circ \chi_{1}^{'}\chi_{0}^{''} \circ
\chi_{0}^{'} \circ  \left(
\sum_{a+b+c=n-3-j}x_{2}^{a}x_{4}^{b}x_{1}^{c})\right)
\left(\sum_{k=0}^{i} x_{1}^{k}x_2^{i-k} \right)\iota'\\

& = & \overline{S} \oplus  \displaystyle\sum_{j=0}^{n-3}
\sum_{i=0}^{n-2} \underbrace{\varepsilon' (x_{1}^{2} - x_{1}x_{2} - x_{1}x_{4} + x_{2}x_{4}) \left(\sum_{a+b+c=n-3-j}x_{2}^{a}x_{4}^{b}x_{1}^{c})\right) \left( \sum_{k=0}^{i} x_{1}^{k}x_2^{i-k} \right)\iota'}_{\Omega}\\
\end{array}$$
\\
where \begin{equation} \overline{S} = S \circ \chi_{0}^{''} \circ
\left( \sum_{i=0}^{n-2} \sum_{k=0}^{i} x_{1}^{k}x_2^{i-k} \chi_{0}^{'}
\circ \iota' \otimes Id \right)
\end{equation}
\\

 To go from line 4 to 5 we recall what these $\chi$'s are:\\
 \begin{figure}[h]
\centerline{
\includegraphics[scale=.7]{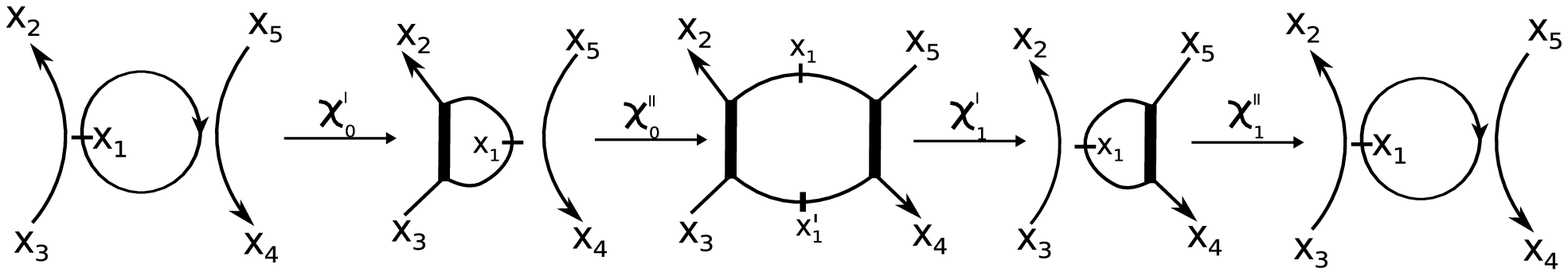}}
\end{figure}

The composition $\chi_{1}^{''} \circ \chi_{0}^{''} \circ \chi_{1}^{'}
\circ \chi_{0}^{'} = (x_{4} - x_{1})(x_{2}-x_{1})=x_{1}^{2} -
x_{1}x_{2} - x_{1}x_{4} +x_{2}x_{4}$, so now we just have to figure
what happens
with $\Omega$.\\

\textbf{Claim} If $i<j$ then $\Omega = 0$ and if $i=j$ then  $\Omega=Id$\\
\indent \emph{Proof:} This is just a simple check. The only thing to
note is that $\Omega
\neq 0$ iff one of the following occurs: \\
$1) c+k=n-1$\\
$2) c+k+1=n-1$\\
$3) c+k+2=n-1$\\
\indent So $i<j \Rightarrow k<j$ so say $c+k=n-1$. Then $a+b+c =
a+b+n-1-k = n-3-j \Rightarrow a+b = -2 + k -j <0$ contradiction,
since $a, b, c$
are nonnegative integers. The other two cases are similar.\\
\indent From above we see that we need $k$ at least equal to $j$. So
if $i=j=k$ and $c+k+2=n-1 \Rightarrow a+b+c=a+b+n-3-k=n-3-j
\Rightarrow a+b=0$ and $\Omega = Id$. The other two cases force
$a+b < 0$. $\square$ \\

So the matrix for $\Omega$ looks like:\\

\centerline {\begin{small}$\left[
\begin{array}{cccclcr}
Id & * & * &  *** &   *** & * & *\\
0  &  Id    & * & *** & *** & * & *\\
\vdots & 0 & \ddots & \ddots &  & \vdots & \vdots\\
\vdots & &  & \ddots & \ddots & \vdots &  \vdots\\
\vdots & & & & \ddots & * & *\\
 0 &  & \ldots & \ldots &  0 & Id & * \end{array} \right]$ \end{small}}

\medskip

Using Gaussian Elimination we see that only the entry
 corresponding to $i=n-2$ survives and the original complex is
 homotopic to:

\begin{figure}[h]
\centerline{
\includegraphics[scale=.6]{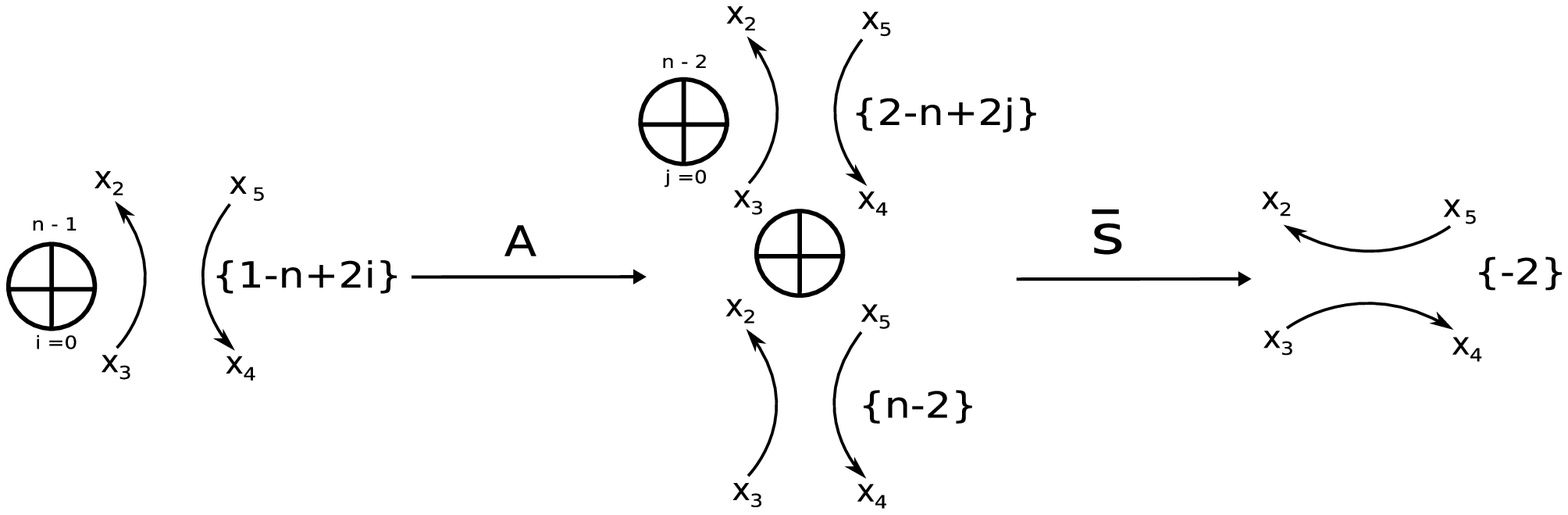}}
\end{figure}

where $A=$

\centerline {\begin{small} $\left[
\begin{array}{cccclcr}
Id & -x_{4} & 0 &  \ldots & \ldots &  0 & 0\\
0  &  Id    & -x_{4} & 0 & \ldots & 0 & 0 \\
\vdots & & \ddots & \ddots & & \vdots & \vdots\\
\vdots & &  & \ddots & \ddots& \vdots & \vdots\\
\vdots & & & & Id & -x_{4}& 0\\
 0 &  & \ldots & \ldots & 0 & Id & -x_{4}\\
 0 &  & \ldots & \ldots & 0 & -Id & x_{2}
 \end{array} \right]$ \end{small}}

\bigskip

This is just our original matrix $\Theta$ but with one more row for
the extra term, for which the entries are computed identically as we
have already done. We reduce the complex in fig. \ref{basictangle}, insert the
overall grading shifts and arrive at our desired conclusion, i.e.:
\medskip
\begin{figure}[h]
\centerline{
\includegraphics[scale=.85]{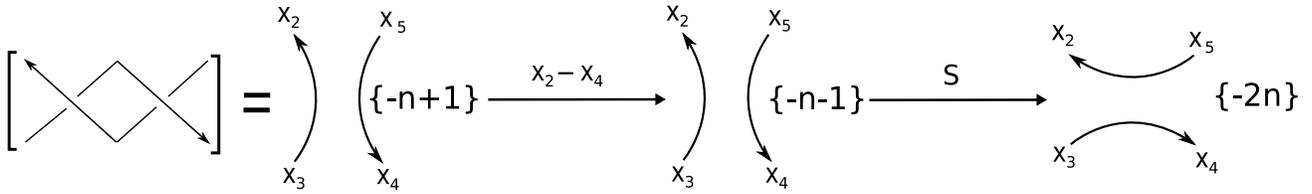}}
\caption{The reduced complex for tangle T}
\label{reducedT}
\end{figure}

Note: to convince ourselves that the map $S$ above is indeed the ``saddle" map as prescribed, we need only to know that the hom-space of degree zero maps between the two right-most diagrams above is $1$-dimensional, in the homotopy category, and then argue that the map is nonzero. This can be done by say closing off the two ends of the tangle above such that we have a non-standard diagram of the unknot and looking at the cohomology of the associated complex.
We leave the details to the reader and refer to \cite{KR} for hom-space calculations.    

\bigskip

\pagebreak

\section{Basic Tensor Product Calculation}

\bigskip

We now consider our tangle T composed with itself, i.e. the tangle gotten by taking two copies of T and gluing the rightmost ends of one to the leftmost of the other. On the complex level this corresponds to taking the tensor product of the complex for T with itself while keeping track of the associated markings.

\begin{figure} [!htbp]
\centerline{
\includegraphics[scale = .9]{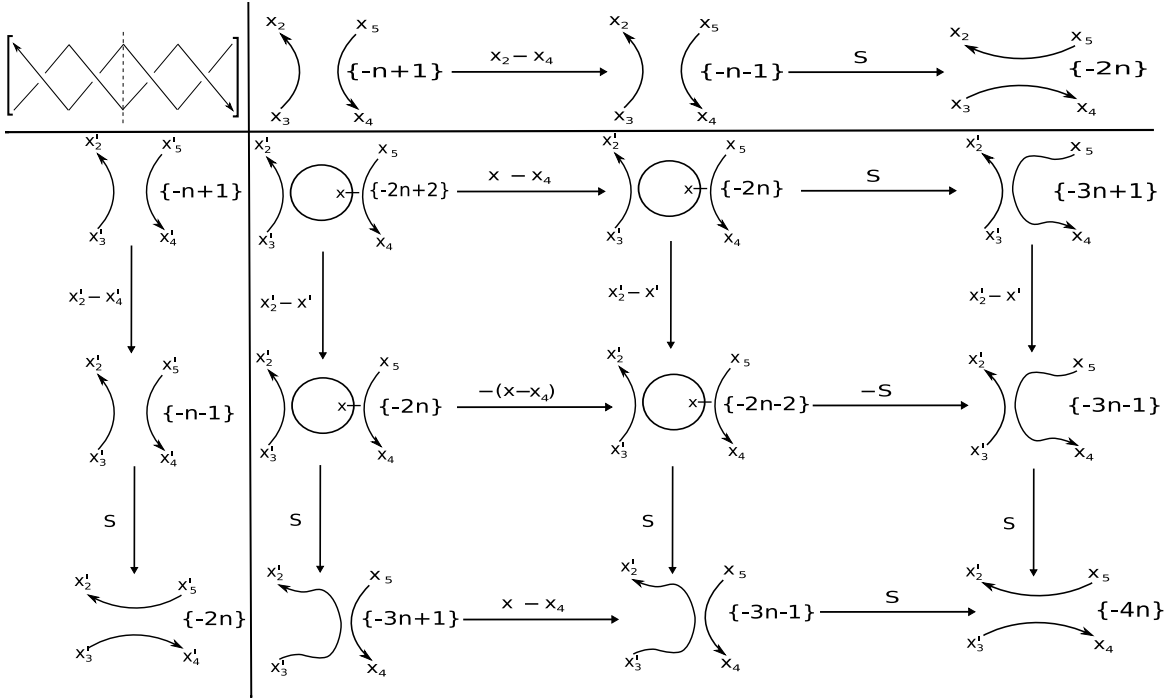}}
\caption{Complex for the tensor product}
\label{tensor}
\end{figure}

\medskip

Note that when we take the tensor product 
we need to keep track of markings. For example:
in the left most entry of the tensored complex $x_2=x'_5=x'_4=x_3$, which
we denote simply by $x$, etc.

As before, we decompose entries in the complex into direct sums of
simpler objects, compute the differentials and reduce using Gaussian
Elimination. In a number of instances we will restrict ourselves to
the $n=3$ case, as the general case works in exactly the same way
with the computation more cumbersome.

We break the computation up based on homological grading.

\pagebreak
 \textbf{Degree 0:}\
\bigskip

\begin{figure} [!htbp]
\centerline{
\includegraphics[scale=.4]{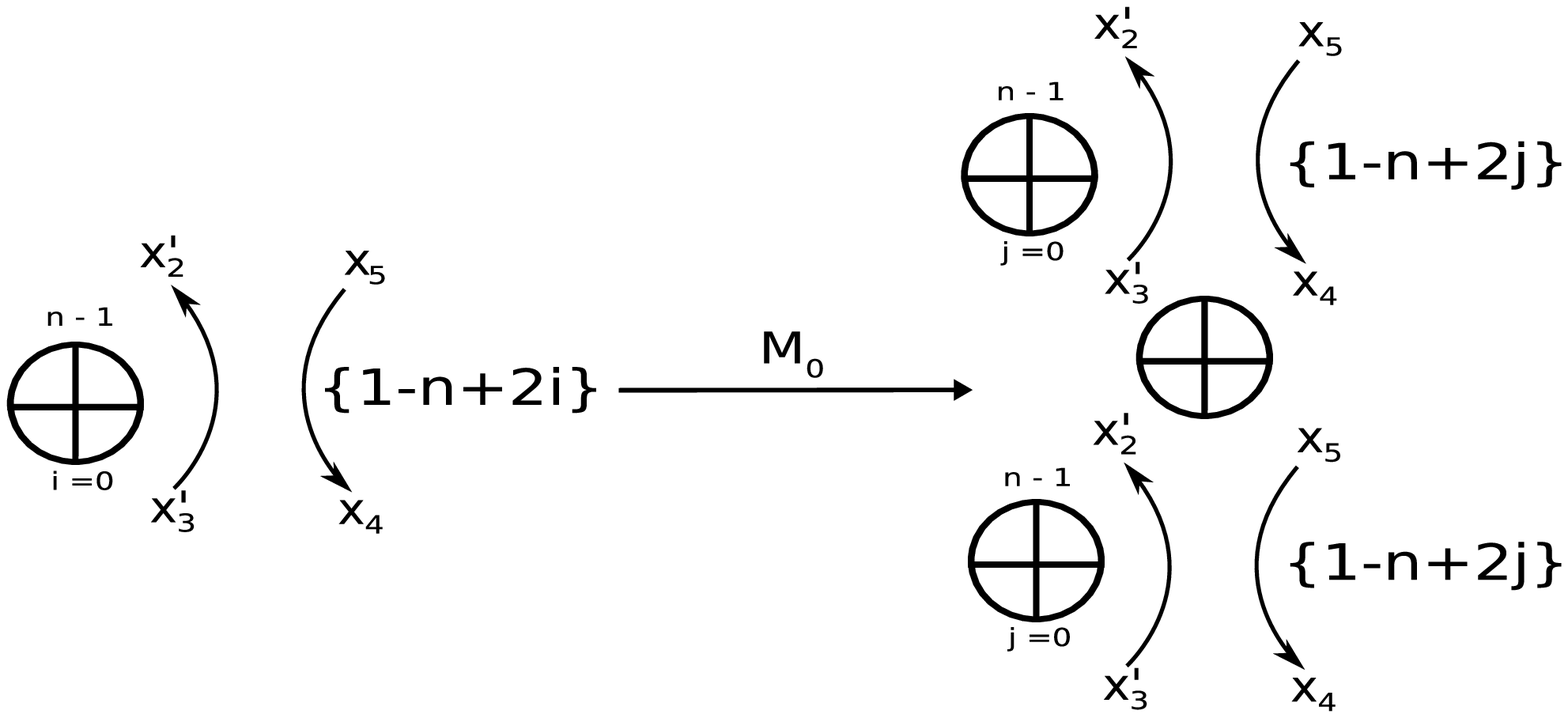}}
\caption{Degree 0 to 1} \label{deg0-1}
\end{figure}\

\bigskip

where $M_{0}$ is:\

\centerline {$\left[
\begin{array}{clcr}
\displaystyle\sum_{i,j=0}^{n-1}Id \otimes \varepsilon(x^{n+i-j} - x^{n-1+i-j}x_{4}) \iota 
\otimes Id\\
\displaystyle\sum_{i,j=0}^{n-1}Id
\otimes\varepsilon(x'_{2}x^{n-1+i-j} - x^{n+i-j})\iota \otimes Id
 \end{array} \right]$}

\medskip

For $n=3$ we have the following:

\medskip

\begin{minipage}[b]{0.2\linewidth}
\centerline{$\left[
\begin{array}{cclcr}
-x_{4} & 0 & 0 \\
Id  &  -x_{4}    & 0\\
0 & Id & -x_{4} \\
x'_{2} & 0 & 0 \\
-Id & x'_{2} & 0 \\
0 & -Id & x'_{2}
\end{array} \right]$}
\end{minipage}
\begin{minipage}[b]{0.1\linewidth}
\begin{center}
$\stackrel{reduce}{\rightsquigarrow}$
\end{center}
\end{minipage}
\begin{minipage}[b]{0.2\linewidth}
\centerline{$\left[
\begin{array}{clcr}
Id &  -x_{4}\\
-x_{4}^{2}    &  0\\
x'_{2}x_{4} & 0\\
x'_{2}-x_{4} & 0\\
-Id & x'_{2}
\end{array} \right]$}
\end{minipage}
\begin{minipage}[b]{0.1\linewidth}
\begin{center}
$\stackrel{reduce}{\rightsquigarrow}$
\end{center}
\end{minipage}
\begin{minipage}[b]{0.2\linewidth}
\centerline{$\left[
\begin{array}{clcr}
0 \\
x'_{2}x_{4}^{2}\\
x'_{2}x_{4}-x_{4}^{2}\\
x'_{2}-x_{4}
\end{array} \right]$}
\end{minipage}
\begin{minipage}[b]{0.1\linewidth}
$= \overline{M}_{0}$
\end{minipage}

\medskip
[Note: we first permute the rows in the first half of the matrix
s.t. the Id maps appear on the diagonal.]

\medskip

The general case is exactly the same, i.e. in the left most matrix above, the upper and lower $3 \times 3$ matrices become expanded to similar $n \times n$ matrices. Hence, the complex reduces to:

\begin{small}
\begin{figure} [!htbp]
\centerline{
\includegraphics[scale=.7]{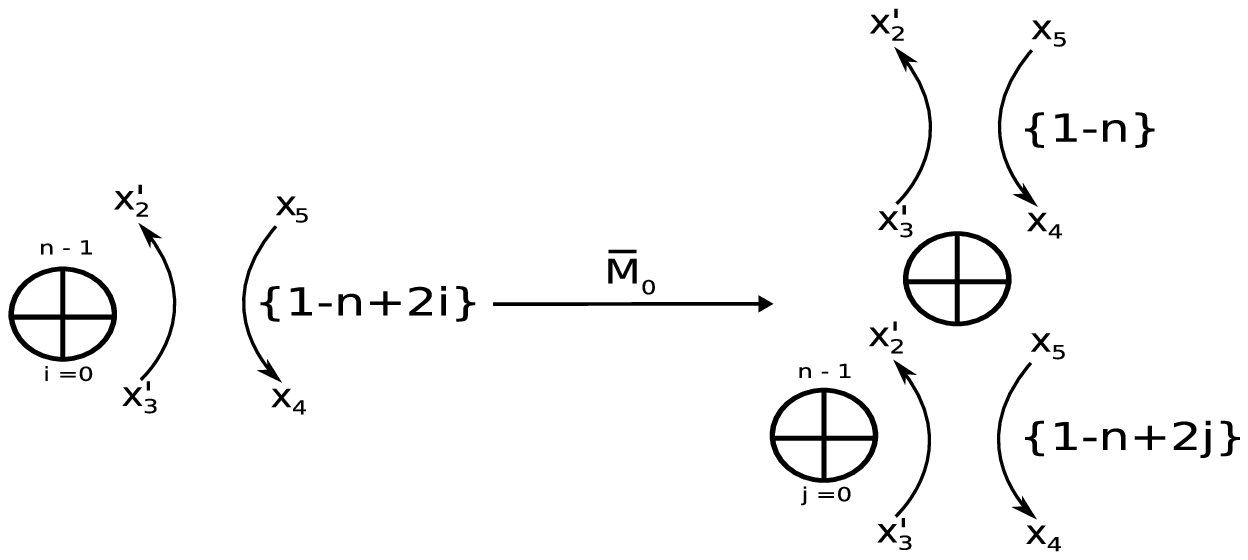}}
\caption{Degree 0 to 1} \label{deg0-1b}
\end{figure}\
\end{small}
\bigskip

\pagebreak
 \textbf{Degree 1:}\

\begin{figure} [!htbp]
\centerline{
\includegraphics[scale = .8]{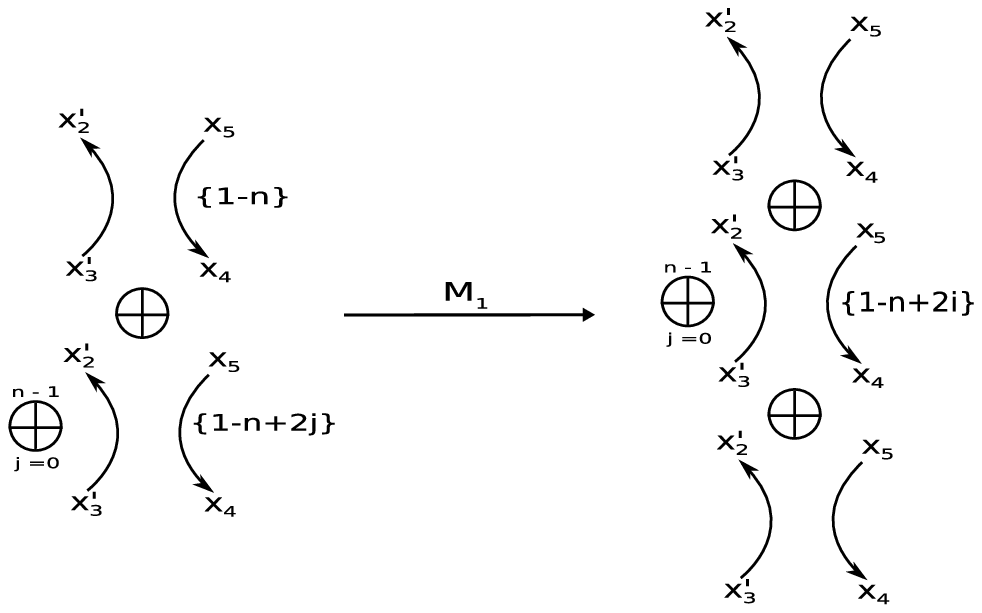}}
\caption{Degree 1 to 2} \label{deg1-2}
\end{figure}\

\medskip

 with $M_{1}$=:\

\medskip

\centerline{\begin{small}$\left[
\begin{array}{cclcr}
Id \otimes S \circ \iota \otimes Id & \{0\}_{1 \times n} \\
M_{1}^{a} & M_{1}^{b} \\
\{0\}_{n \times 1} & M_{1}^{c}
\end{array} \right]$\end{small}}

\medskip

where

$$M_{1}^{a} = \displaystyle\sum_{j=0}^{n-1}Id \otimes \varepsilon
(x'_{2}x^{n-1-j} - x^{n-j}) \iota \otimes Id,$$

$$M_{1}^{b} = \displaystyle\sum_{i,j=0}^{n-1} Id
\otimes\varepsilon (x_{4}x^{n-1-j+i} - x^{n-j+i})\iota \otimes Id,$$

$$M_{1}^{c} = \displaystyle\sum_{i=0}^{n-1}Id \otimes x^{i} S \circ \iota
\otimes Id.$$

\medskip
(Note: $x^{i} S \circ \iota$ here is equal to multiplication by $x_{2}'^{i}$) 
expanding we get:

\bigskip
\begin{small}
\begin{minipage}[b]{0.3\linewidth}
\centerline
 {$\left[
\begin{array}{cccclcr}
Id & 0 & \ldots &  \ldots & \ldots &  0\\
x'_{2}  &  x_{4}   & 0 & \ldots  & \ldots & 0\\
-Id & -Id & x_{4} & 0 &  \ldots & \vdots \\
0 & \ldots & \ddots & \ddots & \ldots & \vdots\\
\vdots & \vdots & \vdots & \ddots & \ddots & 0\\
 0 & 0 & \ldots  & 0 & -Id & x_{4}\\
 0 & Id & x'_{2} & \ldots & \ldots & x'^{n-1}_{2}
\end{array}\right]$}
\end{minipage}
\hspace{.07\linewidth}
\begin{minipage}[b]{0.05\linewidth}
\begin{center}
$\stackrel{reduce}{\rightsquigarrow}$
\end{center}
\end{minipage}
\hspace{.05\linewidth}
\begin{minipage}[b]{0.3\linewidth}
\centerline {$\left[
\begin{array}{cccclcr}
x_{4}   & 0 & \ldots  & \ldots & 0\\
-Id & x_{4} & 0 &  \ldots & \vdots \\
0 & \ddots & \ddots & \ldots & \vdots\\
\vdots & \vdots & \ddots & \ddots & 0\\
0 & \ldots  & 0 & -Id & x_{4}\\
Id & x'_{2} & \ldots & \ldots & x'^{n-1}_{2}
\end{array} \right]$}
\end{minipage}
\hspace{.05\linewidth}
\begin{minipage}[b]{0.05\linewidth}
\begin{center}
$\stackrel{row-moves}{\rightsquigarrow}$
\end{center}
\end{minipage}
\end{small}
\medskip

\begin{small}
\begin{minipage}[b]{0.25\linewidth}
\centerline {$\left[
\begin{array}{cccclcr}
-Id & x_{4} & 0 &  \ldots & \vdots \\
0 & -Id & x_{4}& \ldots & \vdots\\
\vdots & \vdots & \ddots & \ddots & 0\\
0 & \ldots  & 0 & -Id & x_{4}\\
x_{4}   & 0 & \ldots  & \ldots & 0\\
Id & x'_{2} & \ldots & \ldots & x'^{n-1}_{2}
\end{array} \right]$}
\end{minipage}
\hspace{.05\linewidth}
\begin{minipage}[b]{0.05\linewidth}
\begin{center}
$\stackrel{reduce}{\rightsquigarrow}$
\end{center}
\end{minipage}
\hspace{.07\linewidth}
\begin{minipage}[b]{0.25\linewidth}
\centerline {$\left[
\begin{array}{ccclcr}
-Id & x_{4}& \ldots & \vdots\\
\vdots & \ddots & \ddots & 0\\
 0 & \ldots   & -Id & x_{4}\\
x_{4}^{2}   &  \ldots  & \ldots & 0\\
(x'_{2}+x_{4}) & \ldots & \ldots & x'^{n-1}_{2}
\end{array} \right]$}
\end{minipage}
\hspace{.07\linewidth}
\begin{minipage}[b]{0.05\linewidth}
\begin{center}
$\stackrel{reduce}{\rightsquigarrow}$
\end{center}
\end{minipage}
\begin{minipage}[b]{0.25\linewidth}
\centerline {$\left[
\begin{array}{clcr}
0 \\
\displaystyle\sum_{i=0}^{n-1}x'^{i}_{2}x_{4}^{n-1-i}
\end{array} \right]$}
\end{minipage}
\end{small}
\medskip

and we have the following:

\begin{figure} [!htbp]
\centerline{
\includegraphics[scale = .7]{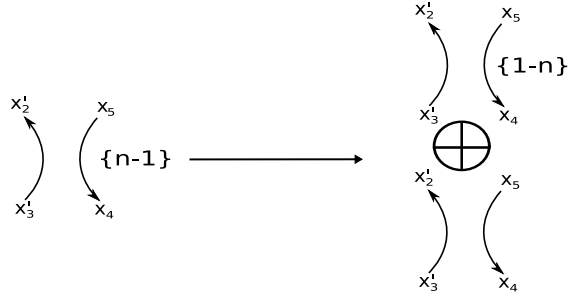}}
\caption{Degree 1 to 2} \label{deg1-2b}
\end{figure}\
\bigskip

 \textbf{Degree 2 and 3:}\

\bigskip

The complex now is pretty simple:
\medskip

\begin{figure} [!htbp]
\centerline{
\includegraphics[scale = .6]{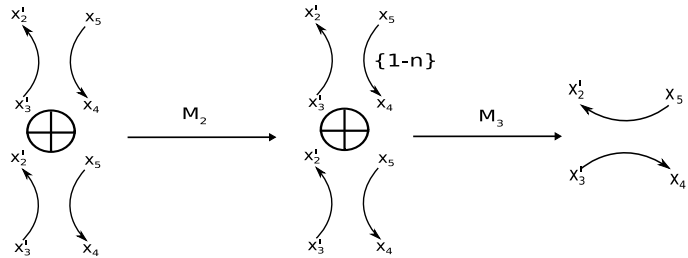}}
\caption{Degree 2 and 3} \label{deg2-3}
\end{figure}\
\medskip

\begin{minipage}[b]{0.05\linewidth}
\begin{center}
$M_{2}=$
\end{center}
\end{minipage}
\hspace{.01\linewidth}
\begin{minipage}[b]{0.5\linewidth}
\centerline {$\left[
\begin{array}{cclcr}
-(Id\otimes S\circ \iota) \otimes Id & Id\otimes (S\circ \iota \otimes Id) \\
0 &  x'_{2} - x_{4}
\end{array} \right],$}
\end{minipage}
\begin{minipage}[b]{0.4\linewidth}
\begin{center}
$M_{3} = \left[ S \hspace{.02\linewidth} S \right].$
\end{center}
\end{minipage}

\medskip
All we have to do is note that $Id\otimes S \circ \iota \otimes Id = Id$ reduce, insert the grading shifts and arrive at the desired
conclusion, i.e.:

\medskip

\begin{figure} [!htbp]
\centerline{
\includegraphics[scale = .9]{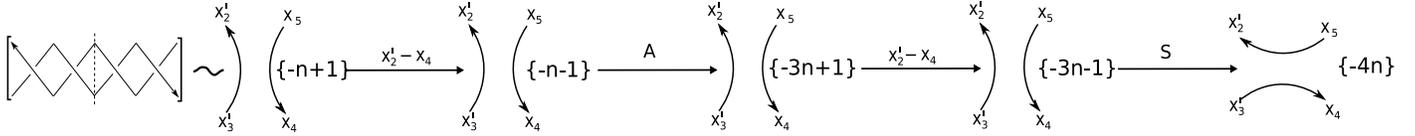}}
\caption{The tensor complex}
\label{basictensor}
\end{figure}\

\medskip

with $A=\displaystyle\sum_{i=0}^{n-1}x_{2}'^{i}x_{4}^{n-1-i}$.

\bigskip
\bigskip

\section{The General Case}
\bigskip

\begin{figure} [!htbp]
\centerline{
\includegraphics[scale = .8]{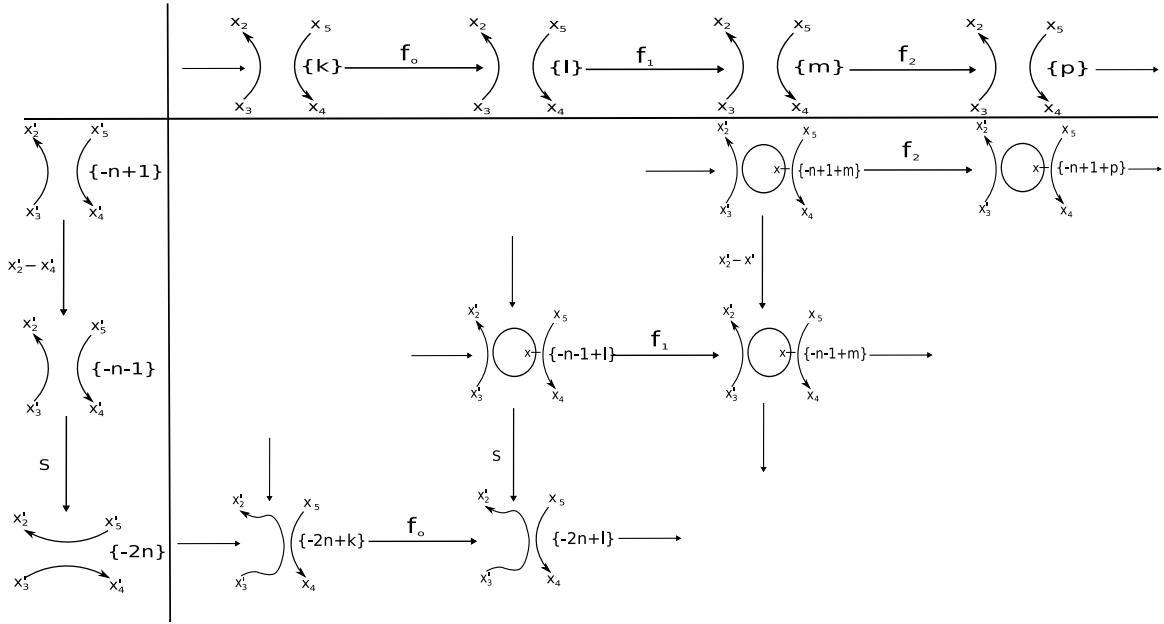}}
\caption{Tensoring the complex with another copy of the basic
tangle $T$}
\label{gentensor} 
\end{figure}\
\bigskip

We suppose by induction that the $k$-fold tensor product of our basic complex has the form as above in fig. $\ref{basictensor}$ with alternating maps $x'_2 - x_4$ and $A$,  the last map being the saddle cobordism $S$, and investigate what happens when
we add one more iteration. As before, this corresponds to tensoring with  another copy of the reduced complex for tangle $T$, i.e. the one in fig. \ref{reducedT}, but as we will see below ``most" of this new complex is null-homotopic and it suffices to consider only the part depicted in fig. \ref{gentensor} directly above.    Note that here the bottom row is a subcomplex which is
isomorphic to that of the top tangle and we claim that, up to
homotopy, this plus two more terms in leftmost homological degree is exactly what survives. The remaining calculation is left to clear up this statement and we begin by taking a look at the highlighted part of the complex depicted in fig. \ref{gentensor}, i.e.:
\pagebreak

\begin{figure} [!htbp]
\centerline{
\includegraphics[scale = .8]{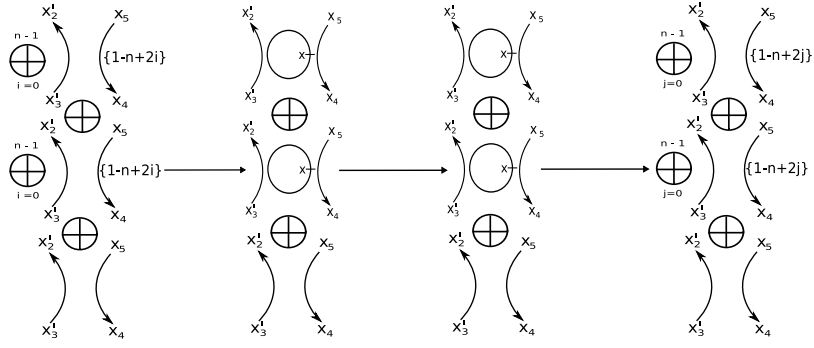}}
\caption{Decomposing the entries of the general tensor product} \label{}
\end{figure}\

...of course we have once again decomposed the complex and left out
the overall grading shifts until later.

The above composition of maps is:

\medskip

\begin{center}{$\left[
\begin{array}{ccclcr}
M^a & \{0\}_{n \times n}  & \{0\}_{n \times 1}\\
M^b & M^c & \{0\}_{n \times 1}\\
\{0\}_{1 \times n} & M^d & f_{0}
\end{array} \right]$}
\end{center}

\begin{minipage}[b]{0.5\linewidth}
$$M^a = \displaystyle\sum_{i,j=0}^{n-1}Id \otimes \varepsilon
f_{2}x^{n-1-j+i}\iota \otimes Id$$
$$M^b = \displaystyle\sum_{i,j=0}^{n-1} Id \otimes\varepsilon
x^{n-1-j}(x'_{2} - x)x^{i}\iota \otimes Id$$
\end{minipage}
\begin{minipage}[b]{0.5\linewidth}
$$M^c = -\displaystyle\sum_{i,j=0}^{n-1} Id
\otimes\varepsilon f_{1}x^{n-1-j+i}\iota \otimes Id$$
$$M^d = \displaystyle\sum_{j=0}^{n-1}Id \otimes  x^{n-1-j} S \circ \iota
\otimes Id$$
\end{minipage}

\medskip
Expanding, with $f_{0}=f_{2} = x-x_{4}$ and $f_{1} =
\displaystyle\sum_{m=0}^{n-1}x^{m}x_{4}^{n-1-m}$ we get the
following submatrices:

\medskip

\begin{minipage}[b]{0.45\linewidth}
\centerline {$M^a = \left[
\begin{array}{cccclcr}
-x_{4}   & 0 & \ldots  & \ldots &  0\\
Id & -x_{4} & 0 & \ldots & \vdots \\
0 & Id & -x_{4}& \ldots & \vdots\\
\vdots & \vdots & \ddots & \ddots & 0\\
0 & \ldots  & 0 & Id & -x_{4}
\end{array} \right]$}
\end{minipage}
\begin{minipage}[b]{0.5\linewidth}
\centerline {$M^b = \left[
\begin{array}{cccclcr}
x'_{2}   & 0 & \ldots  & \ldots &  0\\
-Id & x'_{2} & 0 & \ldots & \vdots \\
0 & -Id & x'_{2}& \ldots & \vdots\\
\vdots & \vdots & \ddots & \ddots & 0\\
0 & \ldots  & 0 & -Id & x'_{2}
\end{array} \right]$}
\end{minipage}

\medskip

\hspace{.08\linewidth}
\begin{minipage}[b]{.7\linewidth}
\centerline {$M^c = - \left[
\begin{array}{cccclcr}
x^{n-1}x_{4}^{n-1} & 0 &  \ldots &  \ldots & 0 \\
\ast & x^{n-1}x_{4}^{n-1} &0& \ldots & 0\\
\vdots & \ldots  & \ddots & \ddots & \vdots\\
\ast & \ldots  & \ast   & x^{n-1}x_{4}^{n-1} & 0 \\
Id   & \ast & \ldots  & \ldots &  x^{n-1}x_{4}^{n-1}
\end{array} \right]$}
\end{minipage}

\bigskip

Now this might look like a mess to reduce, but the thing to notice
is that, in the corresponding summand in our decomposition, the
first matrix above kills off all but the topmost degree terms (with
respect to the decomposition-induced grading shifts), whereas the
$Id$ map found in the left-bottom corner of the second kills off
precisely the topmost degree term. As the maps alternate when we
increase cohomological grading and none of the reductions affect the
bottom row (this is easy to see due to the $0$'s found in the first
row), up to homotopy the bottom row remains altered only by a
grading shift.

As far as the beginning and the end of the complex is concerned we
have already done those computations when we looked at the 2-fold
tensor product. Hence, we arrive at our desired conclusion:
\medskip

\begin{figure} [!htbp]
\centerline{
\includegraphics[scale = .9]{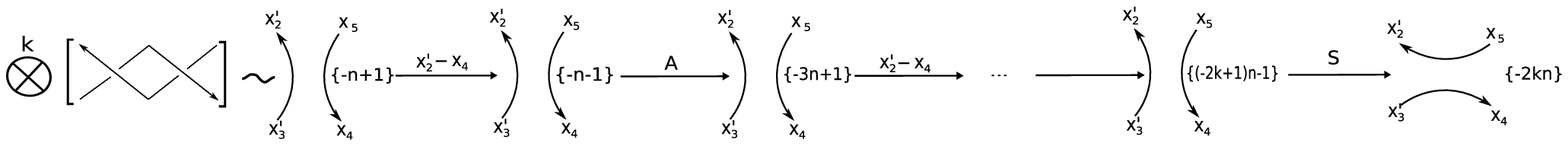}}
\caption{The complex of the k-fold tensor product} \label{}
\end{figure}\
\medskip

where $A = \displaystyle\sum_{i=0}^{n-1} x_{2}'^{i}x_{4}^{n-1-i}$.

Similarly we see that the tangle gotten by flipping all the
crossings is

\begin{figure} [!htbp]
\centerline{
\includegraphics[scale = .9]{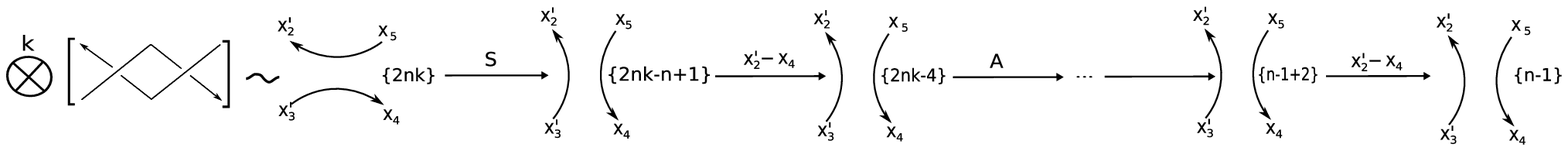}}
\caption{The complex of the k-fold tensor product} 
\end{figure}\
\medskip

\section{Remarks}

Following \cite{BN2} we can propose a similar ``local" algorithm for
computing Khovanov-Rozansky homology. Start with a knot or link
diagram and reduce it locally using the Direct Sum Decompositions
found. Then put all the pieces back together and end up with a
complex where the objects are are just circles, which we can further
reduce to a complex of empty sets with grading shifts, i.e. direct
sums of $\mathbb{Q}$ the maps are matrices with rational entries.
Since a multiplication map $\mathbb{Q} \rightarrow \mathbb{Q}$ is
either an zero or an isomorphism we can use Gaussian Elimination, as
above, to further reduce this complex to one where all the
differentials are zero.  The computational advantage of such an
algorithm is described in more detail in \cite{BN2}. Unfortunately no such program exists to our knowledge. \\
Furthermore, for the examples of tangles we consider here the computational complexity is similar to that of $sl_2$-homology. As there are no more ``thick edges" in any resolution, only Direct Sum Decomposition $0$ is necessary to reduce the complex to $\mathbb{Q}$ vector spaces and matrices between them.  Potentially a modification of the existing programs could allow to compute a large collection of examples composed from these tangles. 

\bigskip

We have done a similar computation for the ``foam" version of
$sl_{3}$-homology introduced in \cite{Kh1}. Here the nodes in the cube of
resolutions are generated by maps from the empty graph to the one at
the corresponding node, with some relations, and the maps are given
by cobordisms between these trivalent graphs. The decompositions
mimic the ones we find here, when specializing to $n=3$, as do the
relations on the maps.  Reducing the complex as before we find that
it is identical to the one found above when specialized to the $n=3$
case. Hence, any link that can be decomposed into the above tangles
has exactly the same homology groups for the "foam" and
matrix-factorization version. This provides a rather vast number of
examples where the isomorphism between the two theories is completely explicit.


\end{document}